\documentclass[oneside,11pt,a4paper]{article}
\usepackage{amsfonts,epsfig,amsmath,amsthm}
\usepackage{verbatim} 

\newtheorem{remark}{Remark}[section]
\newtheorem{definition}{Definition}[section]

\title{On numerical averaging of the conductivity coefficient 
       using two-scale extensions
}

\date{}

\author{Vsevolod Laptev\footnotemark[1]}
\begin{document}
\maketitle
\renewcommand{\thefootnote}{\fnsymbol{footnote}}
\footnotetext[1]{
This work was carried out during the tenures of a
fellowship from University/ITWM in Kaiserslautern (Germany)
and an ERCIM fellowship in Luxembourg and Norway. e-mail: laptevv@mail.ru
}
\begin{abstract}
In this article we compare solutions to elliptic problems having
rapidly oscillated conductivity (permeability, etc) coefficient with
solutions to corresponding homogenized problems obtained from
two-scale extensions of the initial coefficient.
The comparison is done numerically on several one and two dimensional test
problems with randomly generated coefficients for different intensities
of oscillation. The dependency of the approximation error on the size 
of averaging is investigated.
\end{abstract}


%
%%% 1 %%%%%%%%%%%%%%%%%%%%%%%%%%%%%%%%%%%%%%%%%%%%%%%%%%%%%%%%%%%%%%%%%%%%%%%%%
\section{Introduction}
We consider a second order elliptic equation with a 
rapidly oscillated coefficient $a_M(\cdot)$:
\begin{equation}
\label{initial problem}
-\nabla\cdot(a_M(x)\nabla u)=f\qquad\mbox{in $\Omega$,}
\qquad u|_{\partial\Omega}=g.
\end{equation}
The equation appears to describe such problems as the stationary heat transfer in composite 
materials, the flow in non-homogeneous porous media as well as in many others.
For the periodic coefficient, the averaging procedure is well known and
is called the periodic homogenization \cite{BLP},\cite{SP},\cite{JKO}. 
The general non-periodic case is very important for practical
applications (e.g. in geoscience, petroleum engineering), and a vast literature exist
discussing and comparing algorithms intended for averaging the permeability
coefficient (see the reviews 
\cite{Farmer},\cite{RenardMarsily},\cite{WenGomezHernandez},\cite{Effendiev}).
Some algorithms are based on the idea that the effective (averaged, upscaled, equivalent
grid block) permeability field in the whole domain can be determined by solving the flow
problem locally \cite{Durlofsky91}. They vary in the choices of the
local subdomain, the boundary conditions, and the ways to extract the effective permeability 
coefficient from the solution of the local problem. 
These algorithms usually perform well, and intuitively there should
be arguments to justify their usage.

An effective coefficient $A(\cdot)$ in the averaged problem 
\begin{equation}
\label{averaged problem}
-\nabla\cdot(A(x)\nabla U)=f\qquad\mbox{in $\Omega$,}
\qquad U|_{\partial\Omega}=g.
\end{equation}
is generally different from $a_M(\cdot)$, 
and the solution $U$ is different from $u$. 
Therefore there is a trouble with perfect justification:
the difference between these two solutions in some cases can be unacceptable.
Homogenization is known as a rigorous way to justify an averaging process.
This is so because homogenization deals with sequences, not with single problems.
And if the sequence of problems converges in some sense to a limit problem then,
whatever strict requirements we have, there is always a set of problems from the sequence
for which this limit can be considered as an averaged problem.
Although in practice we usually need to upscale a single problem like (\ref{initial problem}), 
not the whole sequence.
Nevertheless, if our initial problem (\ref{initial problem}) belongs to the 
sequence in the
homogenization process then the limit problem may be a reasonable
candidate for upscaled initial problem, even if we cannot improve the
approximation.
We only need that the sequence is homogeneous in the sense that all its members,
including our initial problem, have something in common 
(it is important to avoid situations when a convergent sequence contains an element which has nothing to do with the rest of the sequence). 

One way to do so is to use the sequence from locally periodic
homogenization (see e.g.\cite[p.71]{BLP}) 
\begin{equation}
\label{sequence of problems}
-\nabla\cdot\left(a\left(x,\frac{x}{\varepsilon}\right)\nabla u_\varepsilon\right)=
f\qquad\mbox{in $\Omega$,}
\qquad u_\varepsilon|_{\partial\Omega}=g,
\end{equation}
where the function $a(x,y)$ is a two-scale extension of the initial
coefficient $a_M(\cdot)$:
\begin{definition}[from \cite{tse4npc}]
\label{Def two-scale extension}
Let us say that a function $a(x,y)$, $(x,y)\in\Omega\times\mathbb{R}^d$,
$1$-periodic in the variable $y$, 
is a two--scale extension for $a_M(x)$ if there exists a positive number
$\bar\varepsilon$ such that  
\begin{equation}
 a\left(x,\frac{x}{\bar\varepsilon}\right)
=a_M(x),\qquad \forall x\in\Omega.
\label{two-scale extension}
\end{equation}
\end{definition}
Having a two-scale extension, we can choose a strictly positive
sequence $\{\varepsilon_n\}\to 0$, containing $\bar\varepsilon$, and 
consider (\ref{sequence of problems}) with 
$\varepsilon$ from $\{\varepsilon_n\}$ 
as a sequence in the scope of locally periodic homogenization. 
The expressions for $A(x)$ and corrections of $U$ can be found in the
literature devoted to homogenization. All the members in the sequence
(\ref{sequence of problems}) have in common the function $a(x,y)$;
and at $\varepsilon=\bar\varepsilon$ we recover the initial problem 
(\ref{initial problem}).
In this sense a two-scale extension establishes a connection between
(\ref{averaged problem}) calculated from the homogenization algorithm, and
the initial problem (\ref{initial problem}). As it was already
mentioned, we cannot claim that (\ref{averaged problem}) with such $A(\cdot)$ 
is the averaged problem for (\ref{initial problem}). Moreover, there are
(infinitely) many two-scale extensions leading to different $A(\cdot)$
for the same $a_M(\cdot)$. Nevertheless, we expect that among them there 
could be classes of extensions appropriate for averaging. 
Therefore it is interesting to test numerically the two-scale extensions 
from \cite{tse4npc} on several model problems with non-periodic coefficients.
In each test we calculate both the solution $u$ and the (corrected) solution $U$ 
and verify whether the solutions are close to each other in any sense.
Such numerical evidence could give an idea about the areas of applicability (if any)
of the approach.

The article is organized as follows. In the next section
several ways to construct two--scale extension for arbitrary initial
coefficients $a_M(x)$ are presented. The section 
\ref{Section Averaging using two-scale extension}
contains cell problems and averaging algorithms from the homogenization
theory. Section \ref{Section Numerical results} consists of numerical
results in 1D for ${\mathcal C}$ and ${\mathcal D}_k$ extensions 
(Subsection \ref{1D numerical tests}),
and ${\mathcal D}_2$-extension in 2D (Subsection \ref{2D numerical tests}).

%--------------------------------------------------------------------
\section{Two-scale extensions}
\label{Section Two-scale extensions}
The two-scale extensions (which we numerically investigate in this article) 
and their properties were presented and discussed in \cite{tse4npc}.

The trivial extension
is given by $a(x,y)=a_M(x)$. More useful extensions can be
constructed in the following way (assuming that $a_M(x)$ is known
in a larger domain $\widetilde\Omega\supset\Omega$ in order to avoid
uncertainties close to $\partial\Omega$):
\begin{itemize}
\item 
we choose $\bar\varepsilon>0$ (small in comparison to the typical size
of $\Omega$);
\item
for each $x\in\Omega$ we choose an $\bar\varepsilon$-cube $W_x$
with sides aligned with the coordinate axes, containing $x$: $x\in W_x$.
We also assume that $\widetilde\Omega$ is large enough:  
$W_x\subset\widetilde\Omega$, $\forall x\in\Omega$ 
($W_x$ is a cubic ''Representative Elementary Volume'' around $x$,
$\bar\varepsilon$ is a size of averaging).\\
It is reasonable to distinguish two main choices of $W_x$ 
(${\mathcal C}$ -- continuous, ${\mathcal D}$ -- discrete):
\begin{itemize}
\item[(${\mathcal C}$)]
 $W_x$ is an $\bar\varepsilon$-cube with the center $x$;
\item[(${\mathcal D}$)]
Having a partition 
$\overline{\Omega}=\bigcup_{j=1}^{N_{\mathcal D}}\overline{\Omega_j}$  
($\Omega_i\cap\Omega_j=\emptyset$, $i\ne j$) 
that each $\Omega_j$ has an $\bar\varepsilon$-cube $W_j$
($\Omega_j\subseteq W_j$, $\hat x_j$ is a center of $W_j$) 
then for each $x\in\Omega_j$ we can define $W_x:=W_j$.
\end{itemize}
\end{itemize}
Now we fix $x\in\Omega$ and construct $a(x,\cdot)$:
\begin{enumerate}
\item 
$\tilde a(x,y)=a_M(y)$, $y\in W_x$;
\item 
$\tilde a(x,y)$ is extended $\bar\varepsilon$-periodically in $y$
to the whole space $\mathbb{R}^d$;
\item
$a(x,y)=\tilde a(x,\bar\varepsilon y)$ is the two-scale
extension.
\end{enumerate}
Depending on the choice of $W_x$ we have ${\mathcal C}$-extensions and
${\mathcal D}$-extensions. 
\begin{remark}
The function $a(x,y)$ is still a two-scale extension if we substitute
the item 1. above by one of the following more weak requirements:
\begin{itemize}
\item
$\tilde a(x,y)=a_M(y)$, $y\in O(x)\subseteq W_x$, where $O(x)$ is a
neighbourhood of $x$;
\item
$\tilde a(x,y)=a_M(y)$, $y=x$;
\end{itemize}
and let $\tilde a(x,y)$ to be free in the rest of $W_x$.

This can be used to modify the coefficient near the boundary of $W_x$
e.g. if we want $a(x,y)$ to be continuous.
The second requirement is so weak that it allows to construct any
two-scale extension satisfying Def.\ref{Def two-scale extension} (without saying how
to do it). Probably we should have something like 
$\tilde a(x,y)\approx a_M(y)$, $y\in O(x)$ not to loose the relation
between $a_M(\cdot)$ and $a(\cdot,\cdot)$ completely.
Anyway, we don't consider these possibilities further in this paper.
\end{remark}

The ${\mathcal D}$-extension depends on the choice of $\{\Omega_j\}$.
If we are going to solve (\ref{averaged problem}) using an unstructured
grid then $\{\Omega_j\}$ could be chosen related to that grid. For example,
if we deal with FEM then each $\Omega_j$ could be a union of one or
more finite elements. Here we will test only one kind of a
subdivision of $\Omega$ into $\{\Omega_j\}$, which is more appropriate
for solving (\ref{averaged problem}) on Cartesian grids:
\begin{definition}
Let $k\ge 1$, $\bar\varepsilon>0$ be given. We divide $\mathbb{R}^d$ into
cubes 
$$
\Box_I=\Bigl(i_1h,(i_1+1)h\Bigr)\times\dots\times\Bigl(i_dh,(i_d+1)h\Bigr),\qquad 
h=\bar\varepsilon/k,\quad
I=(i_1,\dots,i_d)\in\mathbb{Z}^d.
$$
We set $\Omega_j=\Box_{I(j)}\cap\Omega$, where $I(j)$ is some numeration of those cubes
which have a non-empty intersection with $\Omega$,  $j=1,\dots,N_{\mathcal D}$. 
$W_j$ is a cube with the side $\bar\varepsilon=kh$, and the center at the same 
point as the center of $\Box_{I(j)}$.
The ${\mathcal D}$-extension constructed this way let us call a ${\mathcal D}_k$-extension.
\end{definition}
%-------------------------------------------------------------
\section{Averaging using two-scale extension}
\label{Section Averaging using two-scale extension}
The sequence (\ref{sequence of problems}) is well investigated in the homogenization theory. 
It is known that the averaged coefficient $A(x)$ at $x\in\Omega$ in the limit problem 
(\ref{averaged problem}) can be calculated via a so-called cell problem. 
Next we remind different formulations of the cell problem applied to the two-scale extensions 
from Section \ref{Section Two-scale extensions}.
Let us fix an arbitrary $x\in\Omega$.
%------------------------------------
\paragraph{Differential form in $\mathbb{R}^d$:}
\begin{equation}
\label{Cell Problem Rd}
\left\{
\begin{array}{l}
-\nabla_y\cdot\Bigl(a(x,y)\bigl(\nabla_y w_j(x,y)+e_j\bigr)\Bigr)=0
\qquad\mbox{in $\mathbb{R}^d$}\\
\int\limits_Y w_j(x,y)\,dy=0,\qquad\mbox{$w_j(x,y)$ is $1$-periodic in $y$}. 
\end{array}
\right.
\end{equation}
%----------------------------------
\paragraph{Differential form in $Y$.}
Since we prefer to solve the problems in a bounded domain, we can
rewrite them in a differential form in a 
cube $Y=(0,1)^d$:
\begin{equation}
\label{Cell Problem Y}
\left\{
\begin{array}{l}
-\nabla_y\cdot\Bigl(a(x,y)\bigl(\nabla_y w_j(x,y)+e_j\bigr)\Bigr)=0
\qquad\mbox{in Y}\\
\mbox{Boundary conditions on $S^0_i$,$S^1_i$ for all $i=1,\dots,d$}:\\
w_j(x,\cdot)|_{S^0_i}=w_j(x,\cdot)|_{S^1_i}\\
 
e_i\cdot\bigl(a(\nabla_y w_j+e_j)\bigr)(x,\cdot)|_{S^0_i}=
e_i\cdot\bigl(a(\nabla_y w_j+e_j)\bigr)(x,\cdot)|_{S^1_i}\\

\int\limits_Y w_j(x,y)\,dy=0
\end{array}
\right.
\end{equation}
where $S^\alpha_i=\{y\in\overline{Y}: y_i=\alpha\}$.
$w_j(x,\cdot)$ is extended periodically in $y$ from $Y$ to
$\mathbb{R}^d$.
%----------------------------
\paragraph{Variational form:}
find $w_j(x,\cdot)\in H^1_{per}(Y)/\mathbb{R}$ such that
\begin{equation}
\label{Variational Cell Problem}
\int_Y \nabla_y\phi(y)^Ta(x,y)\nabla_y w_j(x,y)\, dy=
-\int_Y\nabla_y\phi(y)^Ta(x,y)e_j\, dy\qquad
\forall\phi\in H^1_{per}(Y)/\mathbb{R}.
\end{equation}

The averaged coefficient $A(x)$ can be calculated from the solutions
$w_j$, $j=1,\dots,d$:
\begin{equation}
\label{Averaged coefficient Aij}
A_{ij}(x)=\int_Y e_i^T a(x,y)\Bigl(\nabla_y w_j(x,y)+e_j\Bigr)\,dy.
\end{equation} 
After solving (\ref{averaged problem}), the solution $U$ could be
corrected (see e.g. \cite[p.76]{BLP}):
\begin{equation}
\label{H1-correction}
\widehat{U}(x)=U(x)+\bar\varepsilon\sum_{j=1}^d w_j
\left(x,\frac{x}{\bar\varepsilon}\right)\frac{\partial U}{\partial x_j}(x).
\end{equation} 
Roughly speaking, $\widehat{U}$ approximates $u$ from (\ref{initial problem}) itself, and $U$
approximates the averaged $u$.

Due to $1$-periodicity of $w_j(x,y)$ in $y$, we can substitute $Y$ in 
(\ref{Cell Problem Y})-(\ref{Averaged coefficient Aij})
by any other $1$-cube
$C=(c^1_{min},c^1_{min}+1)\times\dots\times(c^d_{min},c^d_{min}+1)$.
$S^\alpha_i$ we can redefine as 
$\{y\in\overline{C}: y_i=c^i_{min}+\alpha\}$.
For practical purposes it is convenient to take 
$C=Y_x$, where $Y_x=\{y\in\mathbb{R}^d\mid \bar\varepsilon y\in W_x\}$ 
for each fixed $x\in\Omega$. 
Then for $y\in Y_x$ we have $a(x,y)=\tilde a(x,\bar\varepsilon y)=a_M(\bar\varepsilon y)$.
It is also useful for calculating the correction (\ref{H1-correction})
since $x$ being (always) inside $W_x$ implies $x/\bar\varepsilon\in Y_x$.
Thus, we don't need to store $a(\cdot,\cdot)$ as a function of $d\times d$ variables
-- it is possible
to obtain all necessary information directly from $a_M(\cdot)$. 
The averaging method is local:
the averaged coefficient $A(x)$ and $w_j(x,\cdot)$ depend only on the values of 
$a_M(\cdot)$ in $W_x$, a neighbourhood of $x$.

$W_x$ in the ${\mathcal C}$-extension is changing with the point $x$;
the field $A(\cdot)$ is a result of solving the cell problems at all points from $\Omega$.
This is different from the ${\mathcal D}$-extension, where a finite number of cell problems
has to be solved since $W_x$ is the same in $\Omega_j$ ($W_x=W_j$).  
In this case $A(x)$ has a constant value in each $\{\Omega_j\}$.
We note, that the averaged coefficients from both ${\mathcal C}$ and ${\mathcal D}$
extensions coinside at the points $\hat x_j$, the centers of $W_j$.
We also know that $A(x)$ from the
${\mathcal C}$-extension should be continuous (\cite[Prop. 7.1]{tse4npc}).
Therefore from a practical point of view these extensions could be seen as different 
interpretations of the coefficient $A(x)$ known at the finite number of points $\hat x_j$: 
we can treat the data as a continuous or a piecewise constant function.
The continuous data can be interpolated in space between $\hat x_j$.
The interpolation makes possible the numerical averaging with ${\mathcal C}$-extensions. 
Such averaging needs additional care comparing to the averaging from 
${\mathcal D}$-extensions: if the distribution of $\hat x_j$
is not dense enough in $\Omega$, the interpolated field could be significantly 
different from the exact $A(\cdot)$ (e.g. by missing oscillations).  
\begin{remark}
The ${\mathcal D}_k$-extension, (\ref{Cell Problem Y}),
(\ref{Averaged coefficient Aij}) lead to the averaging algorithm for $A(x)$
proposed in \cite[p. 527]{WenDurlofskyEdwards} as one of several alternative 
upscaling procedures using ''border regions''. 
The early algorithm \cite{Durlofsky91} can be obtained from the ${\mathcal D}_1$-extension
(where $W_j=\Omega_j$).
\end{remark}
\begin{remark}
The corrected approximation $\widehat{U}$, calculated from a
${\mathcal D}$-extension via (\ref{H1-correction}), is not continuous. The jumps on
$\overline{\Omega}_n\cap\overline{\Omega}_m$ are expected due to the
abrupt change of the cell solutions $w_j$ when $x$ goes from
$\Omega_n$ to $\Omega_m$. These jumps are more significant for
${\mathcal D}_k$-extensions with smaller $k$ since for large $k$, $W_n$ and
$W_m$ have a large common volume. 
Here it creates no problem since we use only $L^2$, $L^\infty$ norms for the 
comparison $\widehat{U}$ with $u$. 
Although if one is interested in fluxes or $H^1$ approximations then
the correction in the form (\ref{H1-correction}) is probably a 
bad choice. 
\end{remark}
%----------------------------------------------------------------------------------
\section{Numerical results}
\label{Section Numerical results}
Our main purpose in this section is to solve several model problems
(\ref{initial problem}), (\ref{averaged problem}) semi-analytically (if possible) or
numerically and to compare $u$ with $U$ and $\widehat{U}$.
\subsection{Random number generator}
In most of the numerical examples in this article the coefficient
$a_M(\cdot)$ is defined with the help of a random sequence $\{\xi_i\}$.
To generate the sequence $\{\xi_i\}$ of real numbers we read at each
occasion an $i$-th pair $(b_0,b_1)$ of bytes from the file \cite{RandomFile}  
($b_0,b_1\in \mathbb{Z_+}$, $0\le b_0,b_1\le 255$) and calculate
$$
\xi_i=(b_0+b_1 2^8)/(2^{16}-1),\qquad \xi_i\in[0,1],\quad i=1,2,\dots\quad. 
$$
The first five pairs are $(34,178)$,$(52,184)$,$(220,178)$,$(237,13)$,$(19,247)$.
This approach was chosen since it is easy to
reproduce the sequence on different computer platforms.
%----------------------------------------------------
\subsection{1D tests}
\label{1D numerical tests}
\begin{figure}[p]
\centerline{\includegraphics[width=1.25\linewidth]
{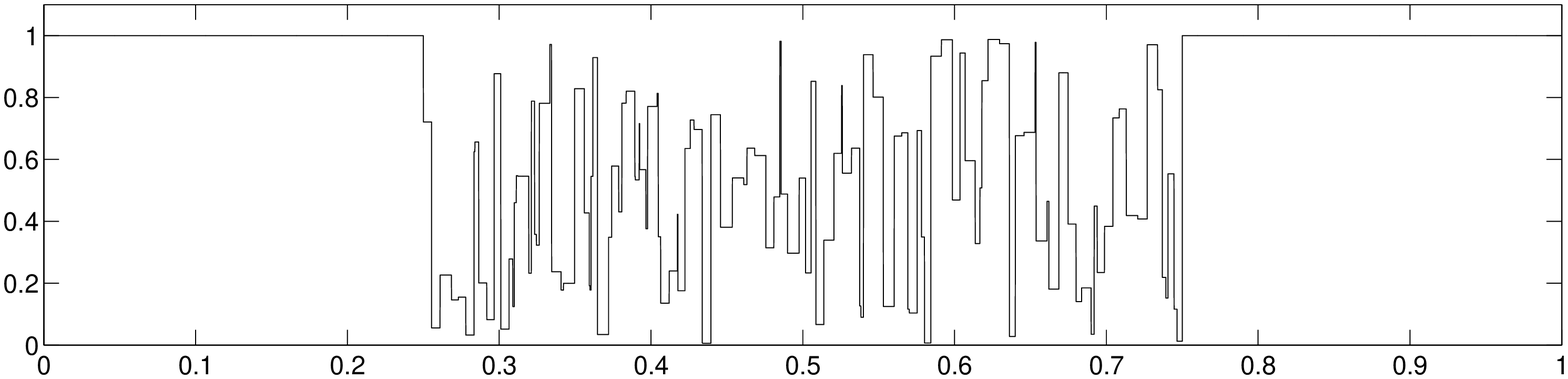}}
\caption{$a_M(\cdot)$ for $\epsilon=0.004$ (case a1)}
\label{D1L125aM}
\end{figure}

\begin{figure}[p]
\centerline{\includegraphics[width=1.25\linewidth]
{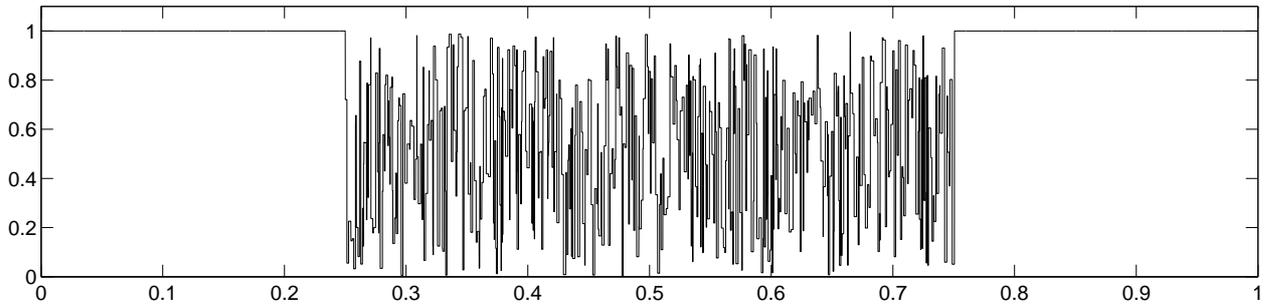}}
\caption{$a_M(\cdot)$ for $\epsilon=0.001$ (case a2)}
\label{D1L500aM}
\end{figure}

\begin{figure}[p]
\centerline{\includegraphics[width=1.25\linewidth]
{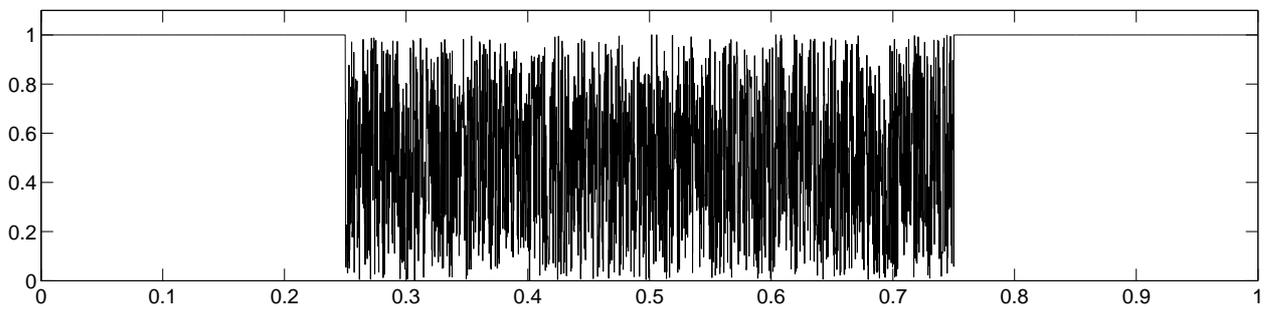}}
\caption{$a_M(\cdot)$ for $\epsilon=0.00025$ (case a3)}
\label{D1L2000aM}
\end{figure}

In 1D we have the following problem
$$
\frac{d}{dx}\left(a(x)\frac{du}{dx}\right)=f(x),\qquad u(0)=u_l,\quad u(1)=u_r,
$$
where $a(x)$ is either the initial coefficient $a_M(x)$ or the averaged 
coefficient $A(x)$. $a_M(x)$ is a constant in $[0,1/4)\cup(3/4,1]$ and has oscillations 
in $[1/4,3/4]$ (see Fig.\ref{D1L125aM}--\ref{D1L2000aM}).
The solution to the equation can be written analytically:
\begin{equation}
\label{exact derivative}
a(x)u'(x)=C+\int_0^xf(x)\,dx=C+F(x),\qquad u'(x)=\frac{C}{a(x)}+\frac{F(x)}{a(x)},
\end{equation}
$$
u(x)=u_l+C\int_0^x\frac{1}{a(x)}\, dx+\int_0^x\frac{F(x)}{a(x)}\, dx,
$$
where $C$ can be determined from the boundary condition $u(1)=u_r$:
$$
C=\left(\int_0^1\frac{1}{a(x)}\, dx\right)^{-1}
\left(u_r-u_l-\int_0^1\frac{F(x)}{a(x)}\, dx\right).
$$
Thus, the semi-analytical numerical solutions for initial and averaged problems need
only the numerical integration.
It seems to be more flexible not to consider $a(\cdot)$ in some exact analytical form,
but to use discretizations of $a(x)$, $u(x)$ on uniform grids.
Thanks to one-dimensionality, grids with millions of
points are available ($N_{sol}$ -- number of points).

The cell problem ($x$ is like a parameter here)
$$
\frac{d}{dy}\left(a(x,y)\Bigl(\frac{dw(x,y)}{dy}+1\Bigr)\right)=0, \qquad
w(x,0)=w(x,1)
$$
also can be solved analytically (up to an additive constant):
$$
\frac{dw(x,y)}{dy}=\frac{C(x)}{a(x,y)}-1,\qquad
w(x,y)=w(x,0)+C(x)\int_0^y\frac{dy}{a(x,y)}-y, 
$$
where $C(x)=\Bigl(\int_0^1a(x,y)^{-1}\, dy\Bigr)^{-1}$, since $w(x,0)=w(x,1)$.
The averaged coefficient is
$$
A(x)=\int_0^1 a(x,y)\left(\frac{dw(x,y)}{dy}+1\right)\, dy=C(x),
$$
$$
A(x)=
\left(\int_0^1 \frac{dy}{\tilde a(x,\bar\varepsilon y)}\right)^{-1}=
\left(\frac{1}{\bar\varepsilon}\int_0^{\bar\varepsilon}\frac{dz}{\tilde a(x,z)}
\right)^{-1}=
\left(\frac{1}{\bar\varepsilon}\int_{w_-(x)}^{w_+(x)}\frac{dz}{a_M(z)}\right)^{-1},
$$
where $W_x=(w_-(x),w_+(x))=
(\hat x(x)-\bar\varepsilon/2,\hat x(x)+\bar\varepsilon/2)$.  
For the ${\mathcal C}$-extension: $\hat x(x)=x$. 
For the ${\mathcal D}_k$-extension: $\hat x(x)=h(\lfloor x/h\rfloor+0.5)$, where $\lfloor y\rfloor$
is the largest number from $\mathbb{Z}$: $\lfloor y\rfloor\le y$.

The $H^1$ correction (\ref{H1-correction}) is
$$
\widehat{U}(x)=U(x)+\bar\varepsilon U'(x)w(x,x/\bar\varepsilon),
$$
where the expression for $U'(x)$ can be found in (\ref{exact derivative}).
\begin{remark}
The harmonic averaging is used in the finite volume method e.g. for
discretizing the elliptic operator with discontinuous coefficients 
\cite{Samarskii}. 
\end{remark}
For the 1D tests the coefficient $a_M(x)$ in $\tilde\Omega=(-1,2)$ is
$$
a_M(x)=\left\{
\begin{array}{cl}
1&x\in(-1,1/4)\cup [x_M,2)\\
0.001+\xi_{2i}& x\in [x_i,x_{i+1}),\quad i=1,\dots,M-1\\
\end{array},
\right.
$$
where
$$
x_1=\frac{1}{4},\quad
x_{i+1}=x_i+\epsilon\frac{(0.1+4\xi_{2i-1})}{2.1},\quad
x_{M-1}<\frac{3}{4}\le x_M,
$$
$\{\xi_i\}$ is the pseudo-random sequence of numbers,
$\epsilon$ is either $0.004$ (case a1, see Fig.\ref{D1L125aM}),
$0.001$ (case a2, see Fig.\ref{D1L500aM}) 
or $0.00025$ (case a3, see Fig.\ref{D1L2000aM}). The homogeneous boundary
conditions $g\equiv 0$ ($u_l=u_r=0$) are chosen. We use three different r.h.s.:
oscillating, constant and discontinuous.
$$
f(x)=\left\{
\begin{array}{cl}
50\sin(30x)&\mbox{case f1}\\
-4&\mbox{case f2}\\
4\bigl(\mathbf{1}_{(1/2,3/4)}-\mathbf{1}_{(1/4,1/2)}\bigr)&\mbox{case f3}\\
\end{array}
\right.,
$$
where $\mathbf{1}_{(a,b)}(x)$ is a characteristic function of $(a,b)$.

\begin{figure}[h]
\centerline{\includegraphics[width=1.25\linewidth]
{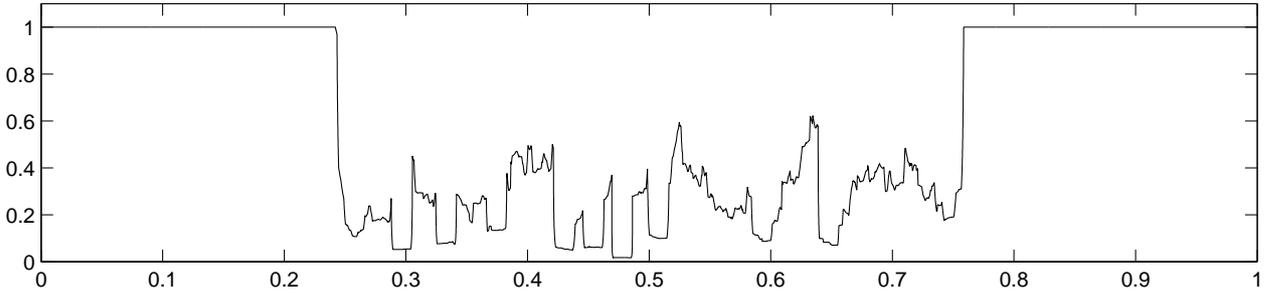}}
\caption{Averaged coefficient $A(\cdot)$ for the case a2 obtained from
${\mathcal C}$-extension for $\bar\varepsilon=0.016$}
\label{D1L500A_Ceps0016}
\end{figure}

\begin{figure}[h]
\centerline{\includegraphics[width=1.25\linewidth]
{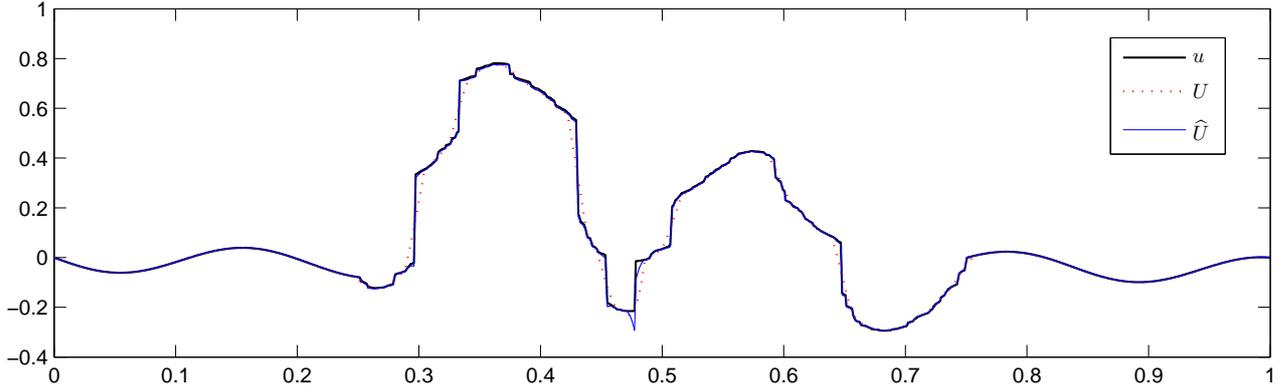}}
\caption{Comparison of $u$ with $U$ and $\widehat{U}$ for
  the case a2f1 obtained from ${\mathcal C}$-extension for $\bar\varepsilon=0.016$}
\label{D1CompL500Ceps0016}
\end{figure}

Let us look at one test more precisely. In Fig.\ref{D1L500A_Ceps0016} the averaged 
coefficient for the case a2, obtained from ${\mathcal C}$-extension for $\bar\varepsilon=0.016$,
is plotted. From Fig.\ref{D1CompL500Ceps0016}, where $U$, $\widehat U$ and $u$
are compared, we see that the averaging is capable to provide good
approximations, and that the correction $\widehat{U}$ approximates $u$ with
a better quality than $U$ (the later looks more like an average of
$u$ smoothing the abrupt curve).

To estimate quantitatively the quality of the approximation we will use:
$$
\widehat{E}_2=\|\widehat{U}-u\|_{L^2(0,1)},\qquad
\widehat{E}_\infty=\|\widehat{U}-u\|_{L^\infty(0,1)},\qquad
E_2=\|U-u\|_{L^2(0,1)},\qquad
E_\infty=\|U-u\|_{L^\infty(0,1)}.
$$
\subsubsection{${\mathcal C}$-extensions in 1D}
\begin{figure}[p]
\centerline{\includegraphics[width=1.25\linewidth]
{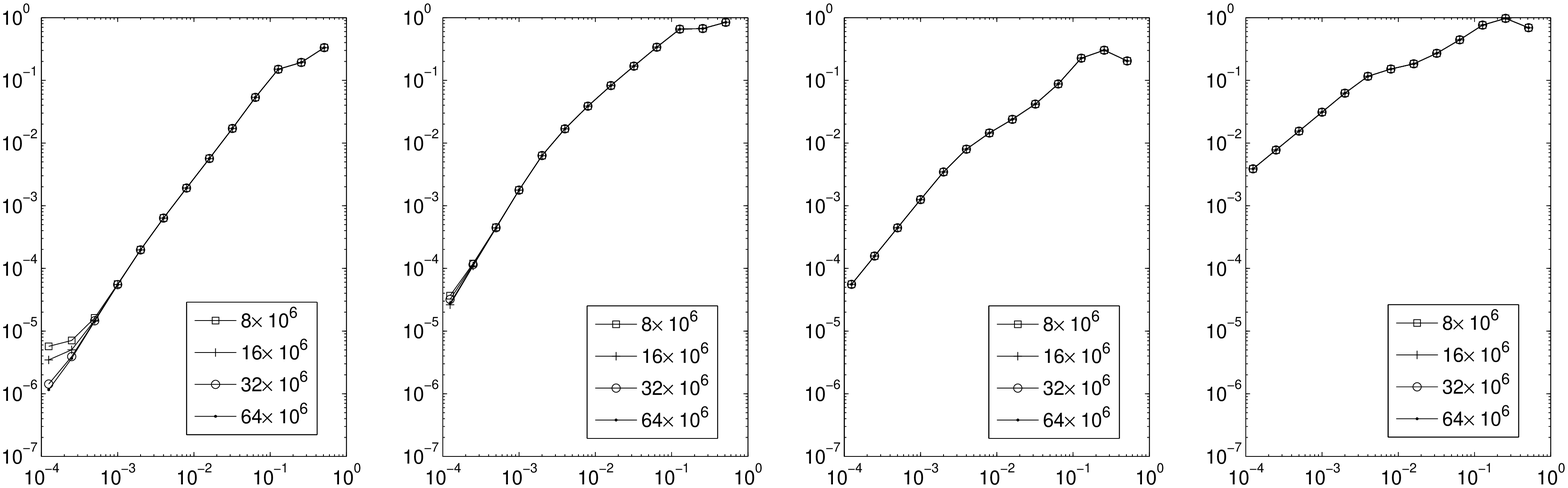}}
\caption{$\epsilon=0.001$, case f1, ${\mathcal C}$-extension: $\widehat{E}_2$, $\widehat{E}_\infty$, $E_2$,
  $E_\infty$ depending on $\bar\varepsilon$}
\label{D1L500F1}
\end{figure}

\begin{figure}[p]
\centerline{\includegraphics[width=1.25\linewidth]
{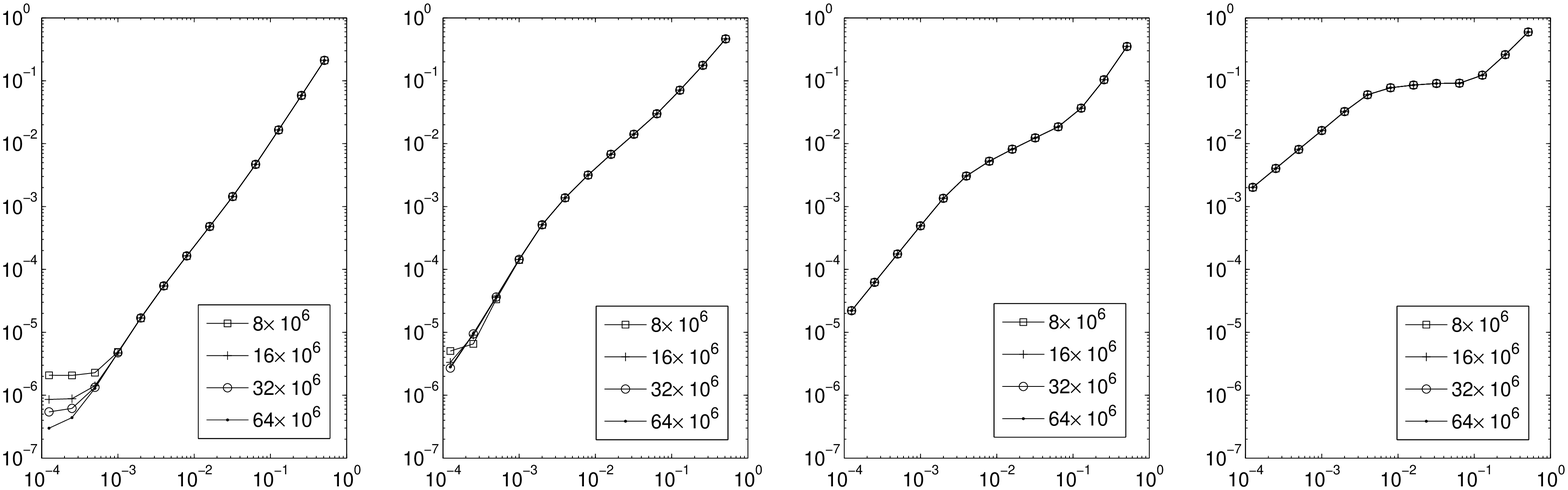}}
\caption{$\epsilon=0.001$, case f2, ${\mathcal C}$-extension: $\widehat{E}_2$, $\widehat{E}_\infty$, $E_2$,
  $E_\infty$ depending on $\bar\varepsilon$}
\label{D1L500F2}
\end{figure}

\begin{figure}[p]
\centerline{\includegraphics[width=1.25\linewidth]
{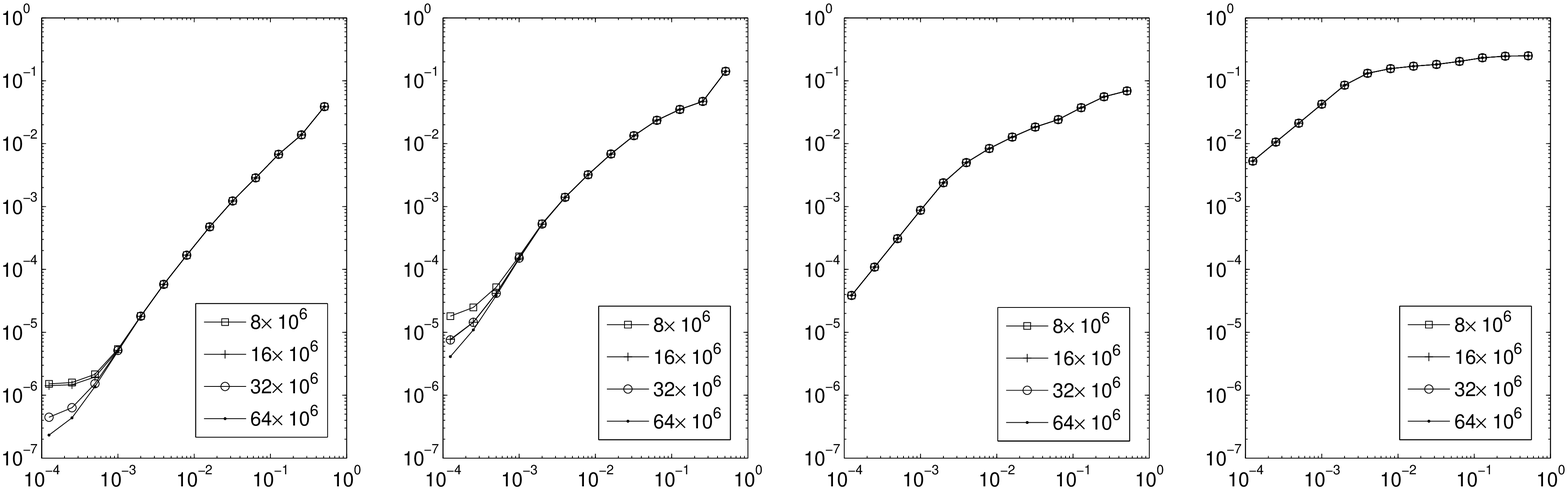}}
\caption{$\epsilon=0.001$, case f3, ${\mathcal C}$-extension: $\widehat{E}_2$, $\widehat{E}_\infty$, $E_2$,
  $E_\infty$ depending on $\bar\varepsilon$}
\label{D1L500F3}
\end{figure}

\begin{figure}[p]
\centerline{\includegraphics[width=1.25\linewidth]
{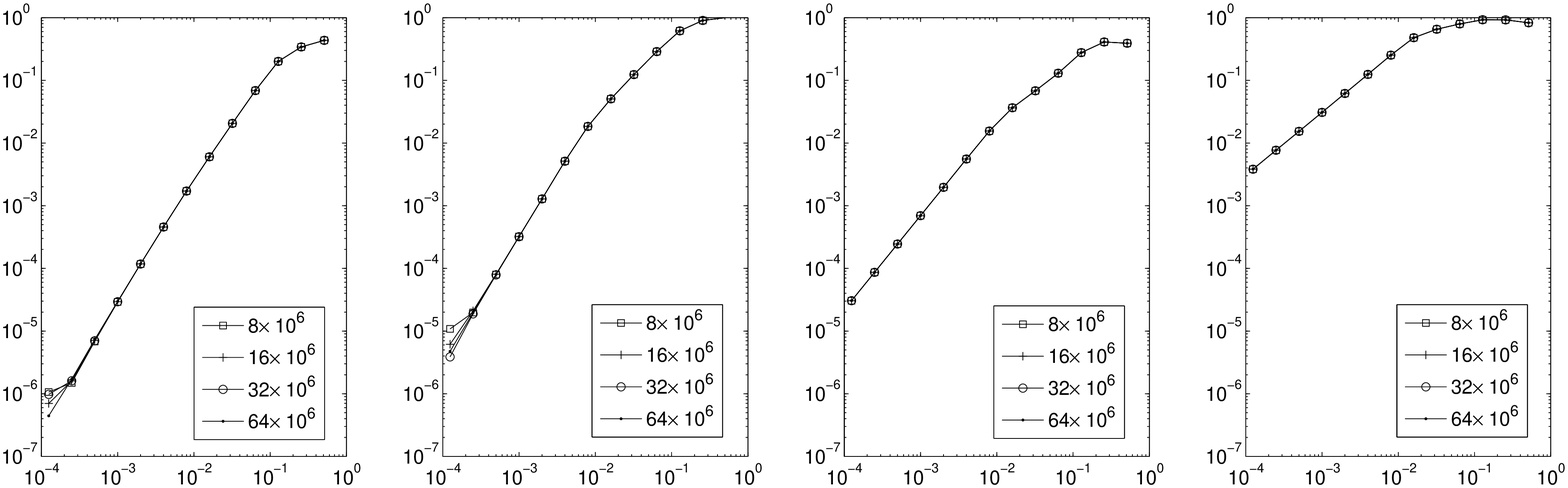}}
\caption{$\epsilon=0.004$, case f1, ${\mathcal C}$-extension: $\widehat{E}_2$, $\widehat{E}_\infty$, $E_2$,
  $E_\infty$ depending on $\bar\varepsilon$}
\label{D1L125F1}
\end{figure}

\begin{figure}[p]
\centerline{\includegraphics[width=1.25\linewidth]
{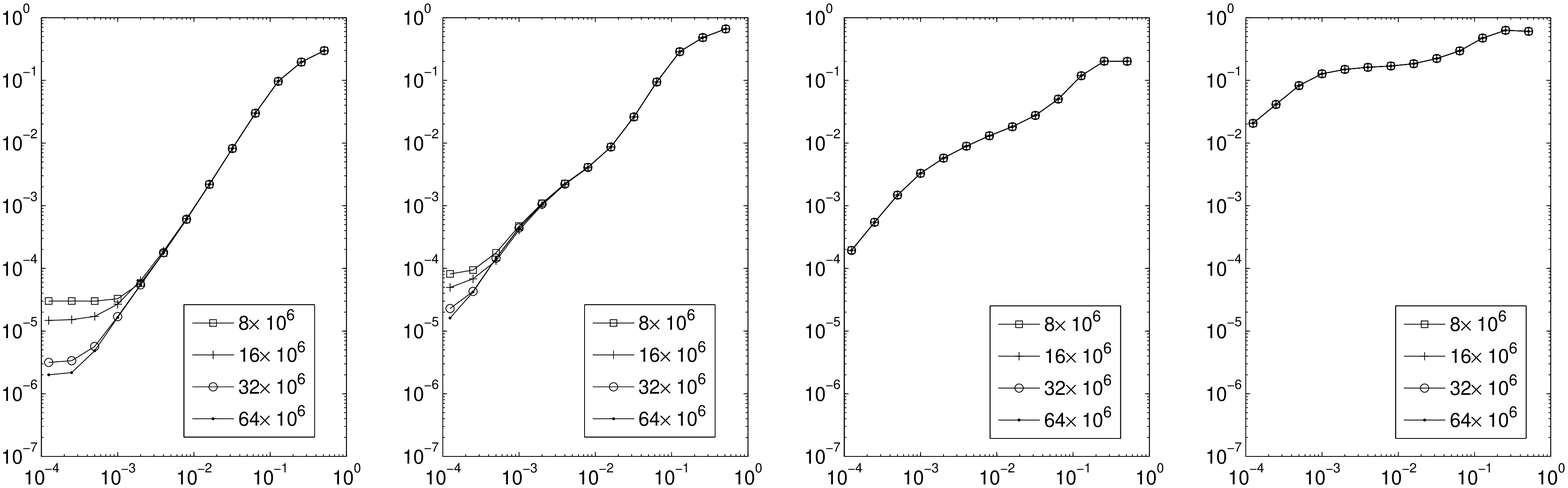}}
\caption{$\epsilon=0.00025$, case f1, ${\mathcal C}$-extension: $\widehat{E}_2$,
  $\widehat{E}_\infty$, $E_2$, 
 $E_\infty$ depending on $\bar\varepsilon$}
\label{D1L2000F1}
\end{figure}

In the first series of tests we solve the problems 
(\ref{initial problem}),(\ref{averaged problem}) for different
$a_M(\cdot)$ and $f(\cdot)$ (cases a2f1, a2f2, a2f3, a1f1, a3f1). The averaged
coefficients are calculated from the ${\mathcal C}$-extensions for different $\bar\varepsilon$. 
The approximation errors are plotted in Fig.\ref{D1L500F1}--\ref{D1L2000F1}.  
In all cases the uniform grids have $N_{sol}=8\cdot 10^6$, 
$16\cdot 10^6$, $32\cdot 10^6$, $64\cdot 10^6$ number of points.
We can see from the figures, that $\widehat{E}_2$, $\widehat{E}_\infty$ 
curves for different $N_{sol}$ are splitted at the end ($\bar\varepsilon\sim 10^{-4}$).
Rounding errors and insufficient resolution could probably explain this, since
the curve obtained on the coarsest grid $N_{sol}=8\cdot 10^6$ starts to
deviate first, and the curve from the finest grid $N_{sol}=64\cdot 10^6$ remains 
longer close to the extrapolated line.
The numerical results show that smaller $\bar\varepsilon$ lead to more accurate approximations,
and that $\widehat U$ approximates $u$ better than $U$ does. The curves on some intervals 
look like straight lines (especially $\widehat{E}_2$). The slopes of the lines on the 
log-log plots give an idea about the order of convergence. 
%---------------------------------example comparison-----------------
%---------------------------------L500-------------------------------
%---------------------------------L125-------------------------------
%---------------------------------L2000-------------------------------
%--------------------------------C vs D_k-------------------------
\begin{figure}[p]
\centerline{\includegraphics[width=1.25\linewidth]
{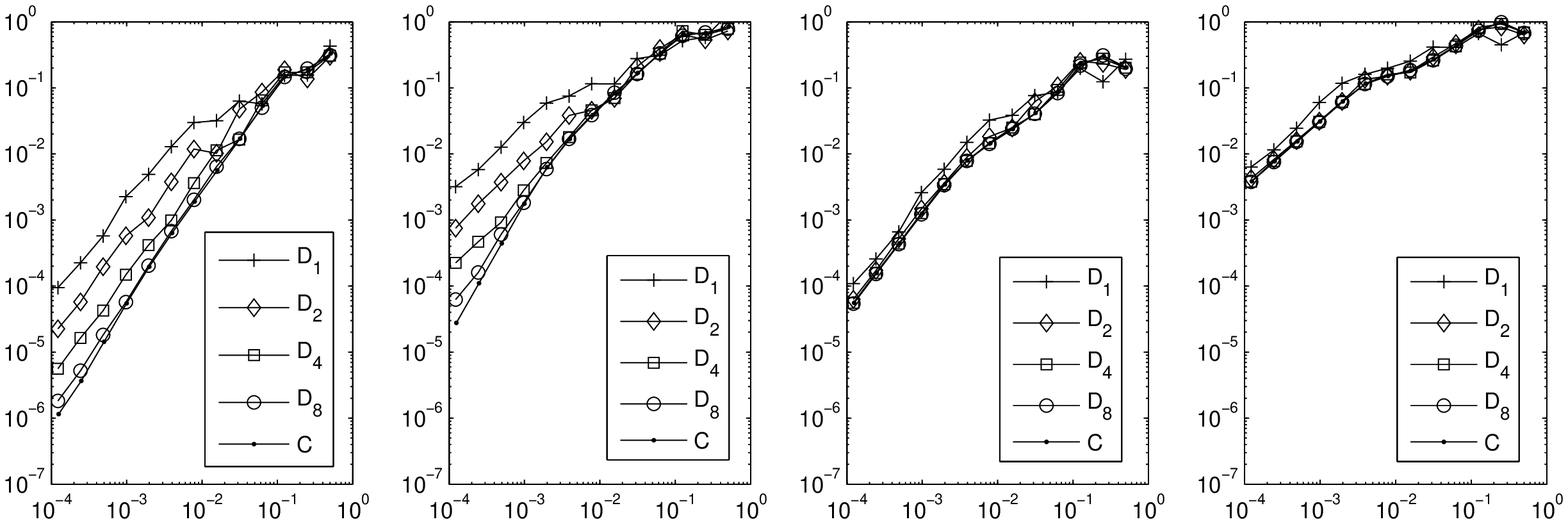}}
\caption{$\epsilon=0.001$, case f1, $\widehat{E}_2$,
  $\widehat{E}_\infty$, $E_2$, $E_\infty$ depending on
  $\bar\varepsilon$ for different extensions:  
${\mathcal D}_1$, ${\mathcal D}_2$,
${\mathcal D}_4$, ${\mathcal D}_8$ and ${\mathcal C}$}
\label{CvsDk}
\end{figure}
%---------------------------------------------------------------------------------
\subsubsection{${\mathcal C}$-extensions and ${\mathcal D}_k$-extensions in 1D}
Calculation of the coefficient $A(\cdot)$ from a ${\mathcal C}$-extension needs high 
computational resources (due to the fine grid), since the fine scale details
of the averaged coefficient (see Fig.\ref{D1L500A_Ceps0016}) could 
disappear after interpolation of a coarse grid data. 
Opposite to that, the averaged coefficient from a ${\mathcal D}$-extension is free from the 
interpolation error, and the needed computational resources are limited by the particular 
choice of the extension. 
Let us compare the qualities of approximation from ${\mathcal C}$, ${\mathcal D}_k$-extensions
for $k=1,2,4,8$. The grid has $N_{sol}=64\cdot 10^6$ nodes. From Fig.\ref{CvsDk} we see
that the ${\mathcal C}$-extension provides better $\widehat U$ approximations (possibly
with higher order of convergence), although there is no significant
difference when $U$ is concerned. We also observe that the quality of approximation from the
${\mathcal D}_k$-extensions approach the quality of approximation from the
${\mathcal C}$-extension when $k$ increases.

The (semi)-analytical solutions $U$, $U'$, $\widehat U$ 
were used also for the ${\mathcal D}_k$-extensions.
This means that the errors which would appear in practical situation ($U_h$,$U_h'$ instead of
$U$,$U'$) were excluded here.
%---------------------------------------------------
\subsection{2D tests}
\label{2D numerical tests}
1D case is very favorable for investigations: extremely fine grids and 
analytical expressions for the solutions are available. 
In 2D we are much more limited in means: we have no analytical
solution for more or less realistic problem specification, and the
finest grid for calculating numerical solutions 
contains only few thousand nodes discretizing OX,OY directions 
(here the maximum is $4096$).
Appearance of arbitrary directions makes the difference from the 1D case.

A reliable investigation of the ${\mathcal C}$-extension remains practically 
out of reach here.  
Thus, we restrict ourselves to ${\mathcal D}_k$-extensions for $k=2$.
The extension has one parameter -- $\bar\varepsilon$. We also use
the equivalent parameter $h=\bar\varepsilon/k=\bar\varepsilon/2$ 
emphasizing that the matrix valued coefficient $A(\cdot)$ is 
a piecewise constant function on the $h$-grid. A coarser grid cannot 
resolve the coefficient properly.

The domain for 2D tests is $\Omega=(0,1)^2$. 
The right hand side and the boundary values for (\ref{initial problem}),
(\ref{averaged problem}) are fixed for all tests: $f(x)\equiv 10$ in
$\Omega$, $g\equiv 0$ on $\partial\Omega$.
The coefficients $a_M(\cdot)$ are described below.
We choose only infinitely smooth coefficients to optimize the accuracy of 
the numerical method on available grids. $a_M(\cdot)$ can be naturally 
extended from $\Omega$ to any $\widetilde{\Omega}\subset\mathbb{R}^2$. 

To solve the 2D elliptic problems with tensor coefficients 
(fine scale problem (\ref{initial problem}), homogenized problem
(\ref{averaged problem}), cell problems (\ref{Variational Cell Problem})) 
we divide the domain $(0,1)^d$ by a uniform Cartesian grid into $N\times N$
squares with the side $h=1/N$ ($h$-grid). 
All squares are subdivided into two triangles by the same diagonal,
and the standard finite element method with linear base functions on such
triangulation is used to solve the problems numerically.
The coefficient is forced to have a constant value inside each square
by taking the value in the center of the square for the whole square
(such approximations are used for (\ref{initial problem}),
(\ref{Variational Cell Problem}) 
since the initial coefficients $a_M(\cdot)$ are smooth in our tests).

The averaged coefficient which is actually used to solve (\ref{averaged problem})
numerically is different from the exact $A(\cdot)$ due to errors
of approximation introduced while solving the cell problems on 
$N_c\times N_c$ grids.
Let us call it 
$A_{h,h_c}(\cdot)$ instead of $A(\cdot)$. The first index 
$h$ emphasizes that the coefficient is piecewise constant on the 
$h$-grid, and the second index $h_c=1/N_c$ specifies the discretization 
step used to solve the cell problems.
$N_c$ is independent from $N=1/h$ and should be large enough for solving 
cell problems with enough accuracy. 
In the tests described below, $N_c$ was usually chosen as large as possible under 
a constrain of reasonable total time of solving $N^2$ cell problems on a 
single processor computer. In addition, the grid 
$(NN_c/k)\times(NN_c/k)$ ($k=2$ here) was fine enough for resolving all
oscillations of $a_M(\cdot)$ in $\Omega$. In some cases $A_{h,h_c}$ 
was compared with $A_{h,2h_c}$, and the solutions of 
(\ref{averaged problem}) with both $A_{h,h_c}$ and $A_{h,2h_c}$ were compared with 
each other in order to verify how the error in $A(\cdot)$ affects the accuracy.  

The problem (\ref{averaged problem}) with the coefficient $A_{h,h_c}(\cdot)$
we solve numerically on two grids: $h$-grid and $h/4$-grid.
The solutions are $U_h$ and $U_{h,4}$ respectively.
$U_h$ is cheap and therefore appropriate for solving practical problems,
although the (coarsest possible) $h$-grid cannot guarantee that $U_h$ 
is a good approximation for $U$. 
For example, the difference between $U_h$ and $U_{h,4}$ is important
when $A_{h,h_c}(\cdot)$ has a high contrast.    
Thus, we need also $U_{h,4}$ -- our numerical substitute for $U$.
 
In order to construct the numerical corrections $\widehat{U}_h$, 
$\widehat{U}_{h,4}$ approximating $\widehat{U}$ from 
(\ref{H1-correction}) we need to save the
solutions of the cell problems. Since the computer memory is also a limited
resource, the cell problem could be solved on $N_c\times N_c$ grid, 
 but saved on $N_{cs}\times N_{cs}$ grid for $N_{cs}\le N_c$.
And we need to store the values of $w_j$ only at the points which 
correspond to $\Omega_i$ inside $W_i$. 
For example, we can choose
a priory a set of points $\{x_k\}$ in $\Omega$ where we would like to know
$\widehat U$, and store the interpolated cell solutions from 
$W_i$ only at the points corresponding to $x_k\in\Omega_i$.
The derivatives from $U$ in (\ref{H1-correction}) are approximated 
in the centers of $h$ squares via central differences and 
then interpolated in $\Omega$. The values in the central differences
are either from $U_h$ or from the projection $U_{h,4}$ to the $h$-grid.

The following relative errors are used to compare the numerical solutions
with the reference solution:
$$
E_2(y)=\|y-u_{ref}\|_{L^2(\Omega)}/\|u_{ref}\|_{L^2(\Omega)},\quad
E_\infty(y)=\|y-u_{ref}\|_{L^\infty(\Omega)}/\|u_{ref}\|_{L^\infty(\Omega)},
$$
where the reference solution $u_{ref}$ is a numerical solution of 
(\ref{initial problem}) obtained on the finest grid $N_{ref}\times N_{ref}$.
$N_{ref}$ is either $2048$ or $4096$ depending on the intensity of 
oscillations in $a_M(\cdot)$.

Each Fig.\ref{Ming Yue test},\ref{log10 E N=64}--\ref{log10 E N=512} 
consists of two subfigures with $E_2$ (left) and $E_\infty$ (right) 
error functions.  
On each subfigure there are 3 functions: $c_1(h)$, $c_2(h)$, $c_3(h)$.
The markers correspond to all test cases.
\begin{itemize}
\item[$c_1$]
The curves with square markers represent the functions
$c_1(h)=E_2(u_h)$ for the left subfigure, and $c_1(h)=E_\infty(u_h)$ for 
the right subfigure, where $u_h$ is the numerical solution of 
(\ref{initial problem}) obtained on the $h$-grid without averaging.  
$u_h$ on the finest grid is the reference solution $u_{ref}$ and therefore
the corresponding square markers for $E_2(u_{ref})=E_\infty(u_{ref})=0$ 
are excluded from the curves.
\item[$c_2$]
The curves with circles represent the functions $c_2(h)=E_2(\widehat U_h)$ 
for the left subfigure, and $c_2(h)=E_\infty(\widehat U_h)$ for the right subfigure. 
\item[$c_3$] 
The curves with point markers represent the functions 
$c_3(h)=E_2(\widehat U_{h,4})$ for the left subfigure, and 
$c_3(h)=E_\infty(\widehat U_{h,4})$ for the right subfigure. 
The averaged coefficient is the same as for $c_2$ -- $A_{h,h_c}(\cdot)$,
but $c_3$ is different from $c_2$. 
\end{itemize}
%-------------------------------------------------------------
\subsubsection{Test with explicitly given coefficient}
\label{subsection MingYue aM}
In \cite{MingYue} the following coefficient for 
(\ref{initial problem}) was proposed as a test ''without scale separation'': 
$$
a_M(x_1,x_2)=\frac{1}{6}\left(
\frac{1.1+\sin(2\pi x_1/\varepsilon_1)}{1.1+\sin(2\pi x_2/\varepsilon_1)}+
\frac{1.1+\sin(2\pi x_2/\varepsilon_2)}{1.1+\cos(2\pi x_1/\varepsilon_2)}+
\frac{1.1+\cos(2\pi x_1/\varepsilon_3)}{1.1+\sin(2\pi x_2/\varepsilon_3)}+
\right.
$$
$$
\left.
+\frac{1.1+\sin(2\pi x_2/\varepsilon_4)}{1.1+\cos(2\pi x_1/\varepsilon_4)}+
\frac{1.1+\cos(2\pi x_1/\varepsilon_5)}{1.1+\sin(2\pi x_2/\varepsilon_5)}+
\sin(4x_1^2x_2^2)+1\right).
$$
where $\varepsilon_1=1/5$, $\varepsilon_2=1/13$, $\varepsilon_3=1/17$, 
$\varepsilon_4=1/31$, $\varepsilon_5=1/65$. 

\begin{figure}[h]
\centerline{\includegraphics[width=1.0\linewidth]
{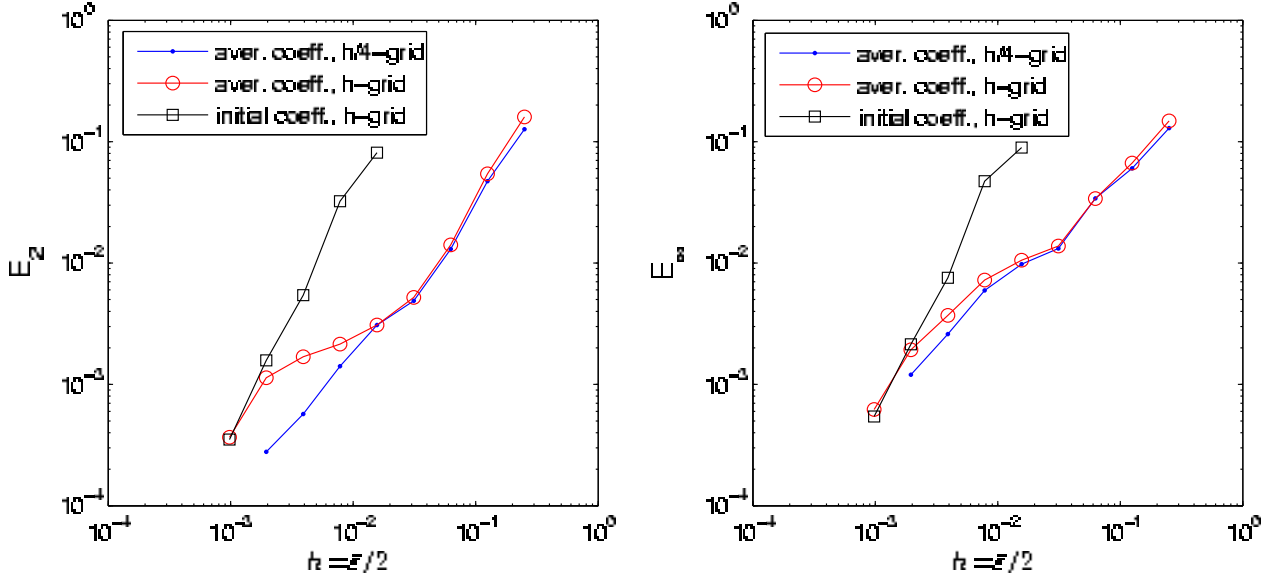}}
\caption{$E_2(\widehat{U}_{h,4})$, $E_2(\widehat{U}_h)$, $E_2(u_h)$ -- left,
$E_\infty(\widehat{U}_{h,4})$, $E_\infty(\widehat{U}_h)$, $E_\infty(u_h)$   
-- right}
\label{Ming Yue test}
\end{figure}

The curves $c_1$,$c_2$,$c_3$ for this test are plotted in 
Fig.\ref{Ming Yue test}.
%------------------------------------------------------------
\subsubsection{Tests with randomly constructed coefficients}
\label{subsection Random aM}
Let us consider the scalar coefficient $a_M(x)=10^{\beta S(x)}$, where
$$
S(x)=\sum_{i=1}^{N_{\sin}}\sin\bigl(\pi i(x_1\sin(\psi_i)+x_2\cos(\psi_i)+\phi_i)\bigr),
\qquad \psi_i=2\pi\xi_{2i-1},\quad
\phi_i=2\xi_{2i}\quad
\beta=\frac{\log_{10}(C)}{M-m},
$$
$\{\xi_i\}$ is the pseudo-random sequence of numbers, the constants $m$,$M$
$$
\begin{array}{ccccc}
N_{\sin}&64&128&256&512\\
m&-19.7229&-36.1412&-49.6262&-81.8554\\
M&22.5351&34.124&51.5507&75.7885
\end{array}
$$
give approximations to minimum and maximum values of $S(x)$ in
$\Omega$ respectively. This allows us to choose the constant $C=10^4$ as the
contrast for $a_M(\cdot)$ ($C\approx\max_x a_M(x)/\min_x a_M(x)$).

\begin{figure}[h]
\centerline{\includegraphics[width=1.2\linewidth]
{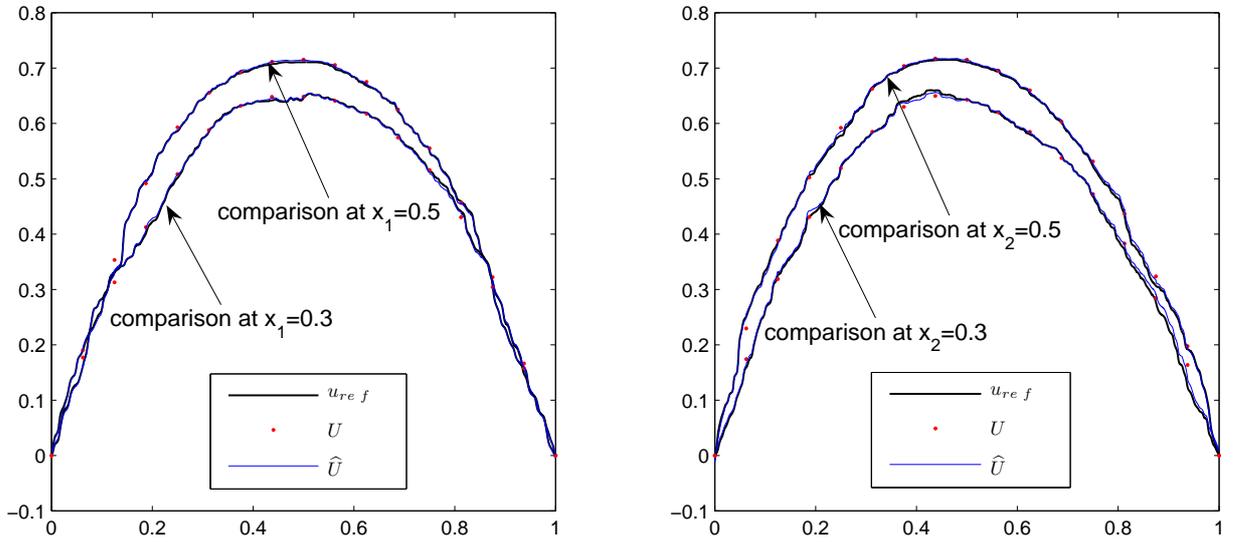}}
\caption{Comparison of $u$ with $U$ and $\widehat{U}$ on several
  cross-sections for $N_{\sin}=256$, $h=1/16$, ${\mathcal
    D}_2$-extension. $u_{ref}$ was calculated on $4096^2$ grid,
  $U$ on $16^2$ grid, cell problems on $512^2$ grids.}
\label{comp N=256 h=1/16}
\end{figure}

\begin{figure}[h]
\centerline{\includegraphics[width=0.9\linewidth]
{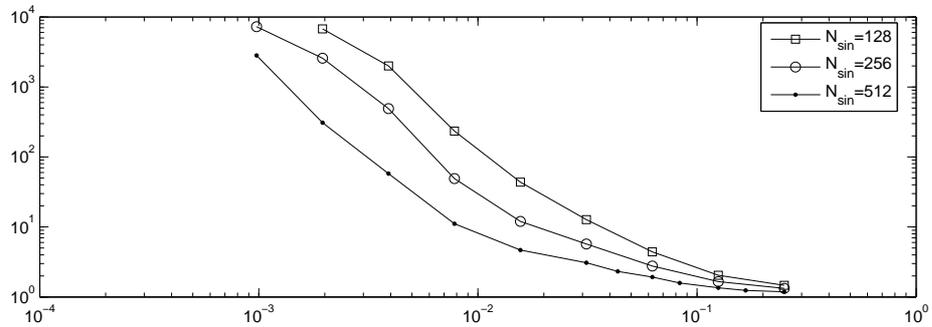}}
\caption{Contrast of the averaged coefficient depending on
  $h=\bar\varepsilon/2$. The contrasts of $a_M$ are $10^4$.}
\label{Contrast2}
\end{figure}

\begin{figure}[p]
\centerline{\includegraphics[width=1.0\linewidth]
{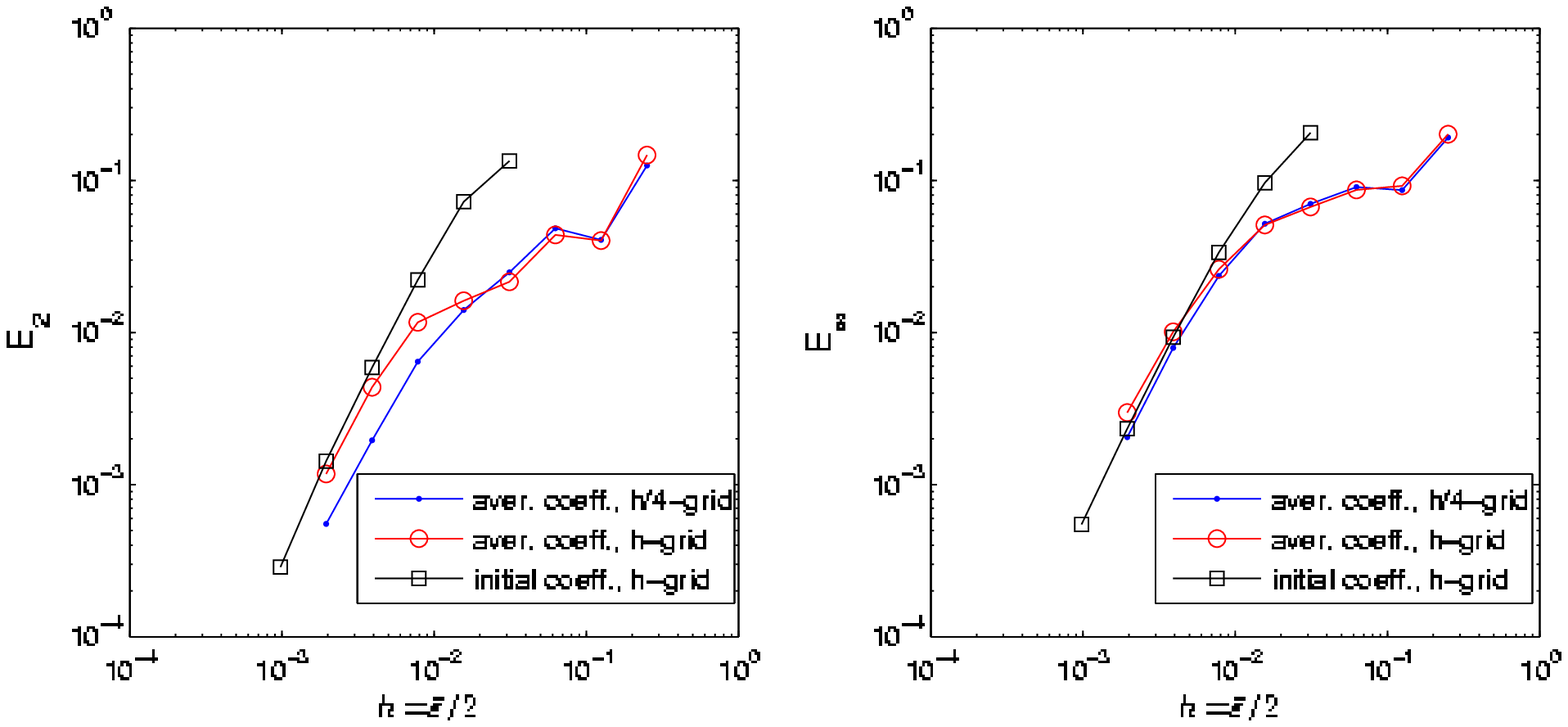}}
\caption{$E_2(\widehat{U}_{h,4})$, $E_2(\widehat{U}_h)$, $E_2(u_h)$ -- left,
$E_\infty(\widehat{U}_{h,4})$, $E_\infty(\widehat{U}_h)$, $E_\infty(u_h)$ 
-- right, $N_{\sin}=64$}
\label{log10 E N=64}
\end{figure}

\begin{figure}[p]
\centerline{\includegraphics[width=1.0\linewidth]
{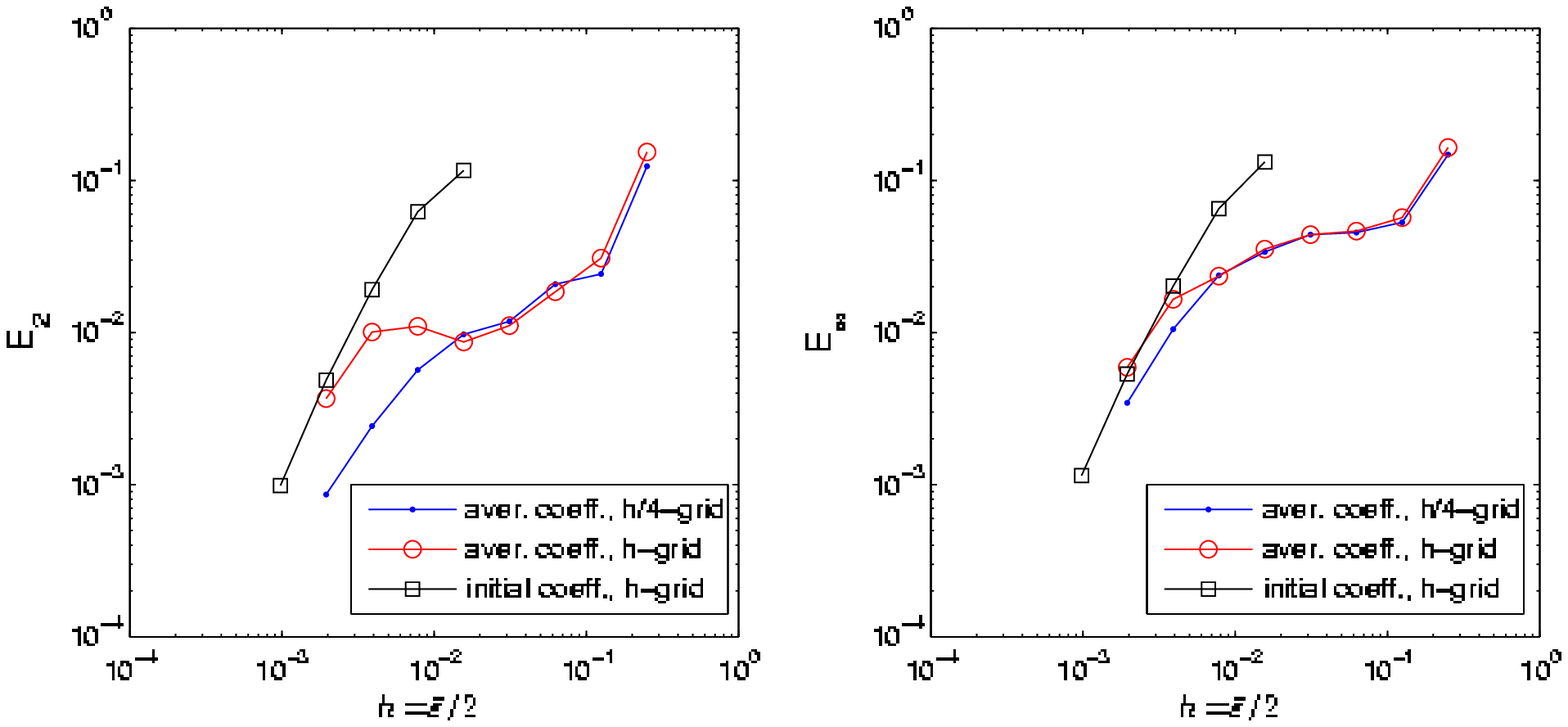}}
\caption{$E_2(\widehat{U}_{h,4})$, $E_2(\widehat{U}_h)$, $E_2(u_h)$ -- left,
$E_\infty(\widehat{U}_{h,4})$, $E_\infty(\widehat{U}_h)$, $E_\infty(u_h)$
-- right, $N_{\sin}=128$}
\label{log10 E N=128}
\end{figure}

\begin{figure}[p]
\centerline{\includegraphics[width=1.0\linewidth]
{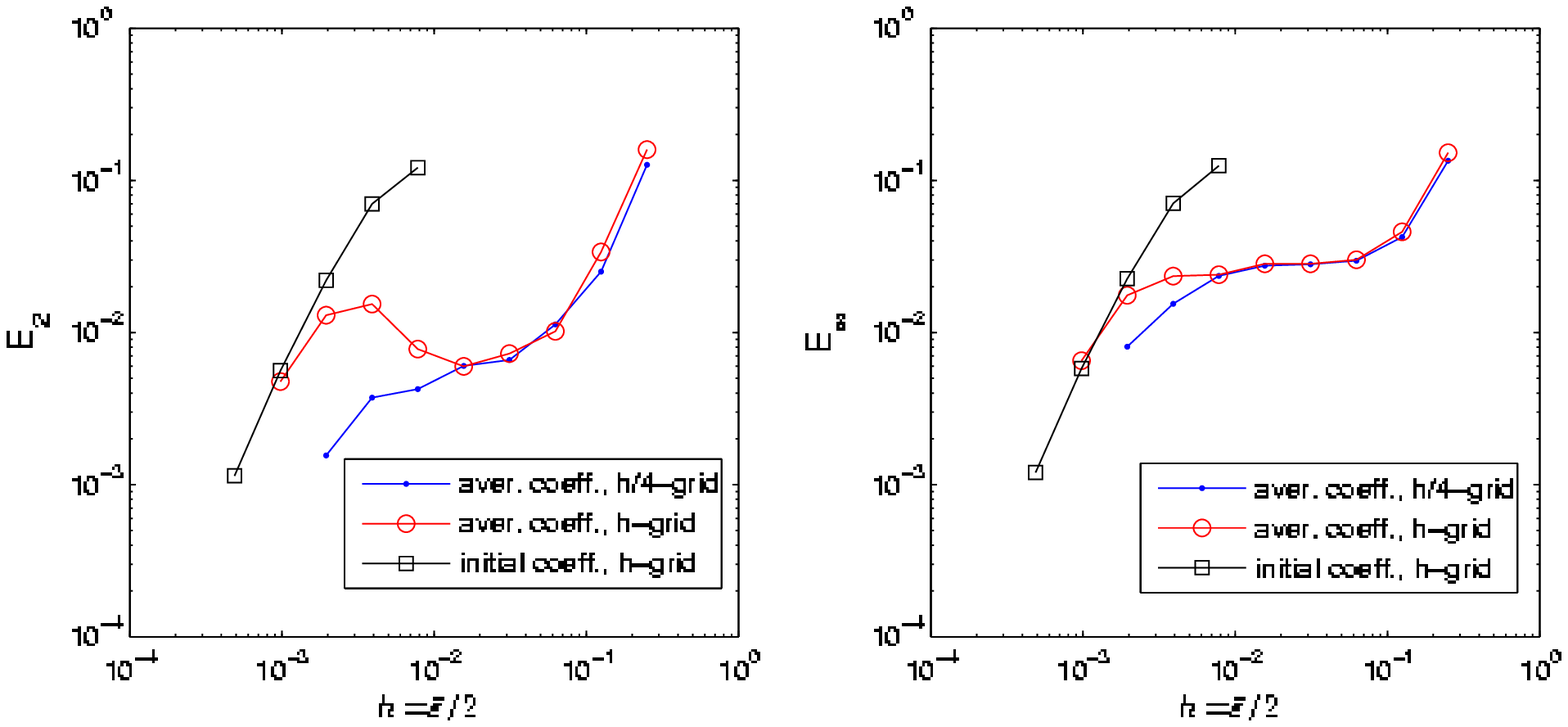}}
\caption{$E_2(\widehat{U}_{h,4})$, $E_2(\widehat{U}_h)$, $E_2(u_h)$ -- left,
$E_\infty(\widehat{U}_{h,4})$, $E_\infty(\widehat{U}_h)$, $E_\infty(u_h)$
 -- right, $N_{\sin}=256$}
\label{log10 E N=256}
\end{figure}

\begin{figure}[p]
\centerline{\includegraphics[width=1.0\linewidth]
{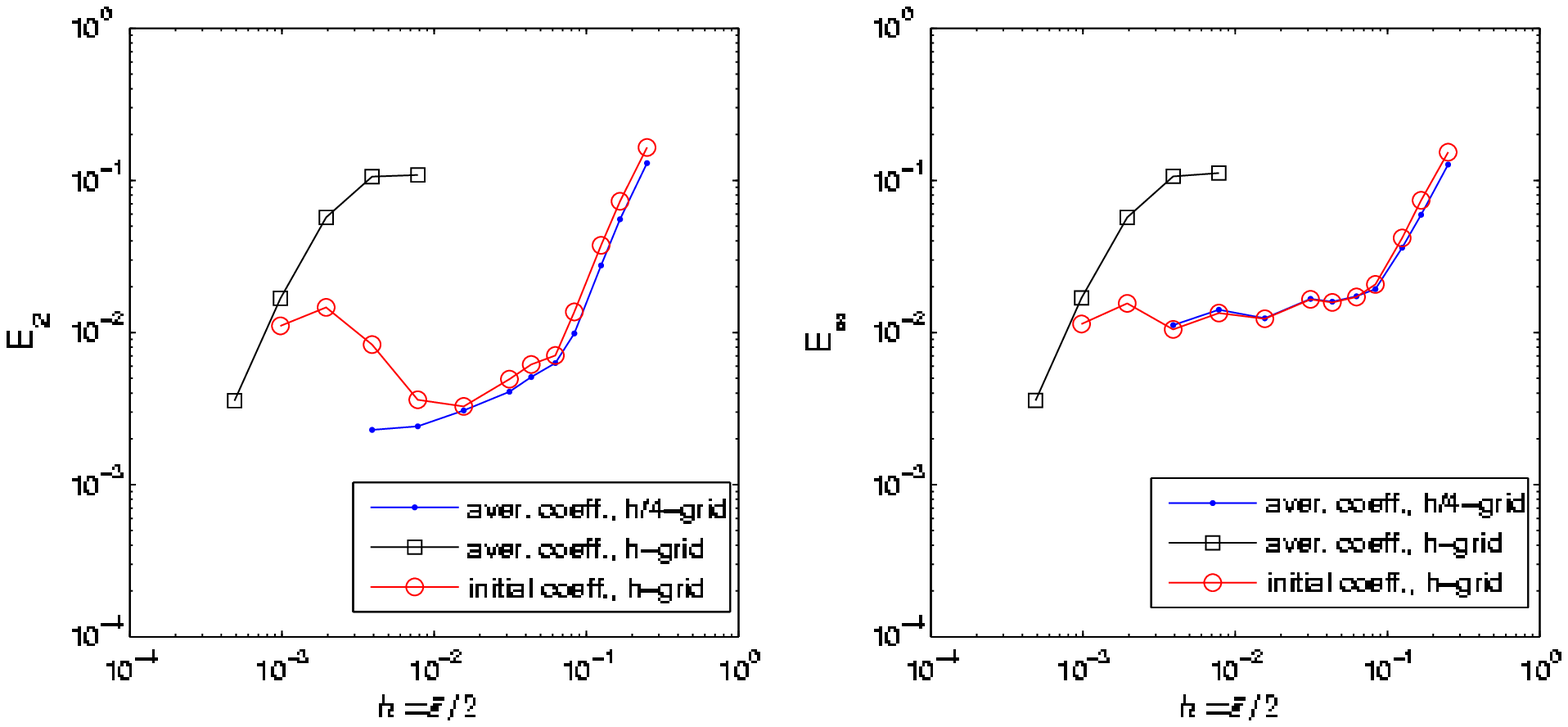}}
\caption{$E_2(\widehat{U}_{h,4})$, $E_2(\widehat{U}_h)$, $E_2(u_h)$ -- left,
$E_\infty(\widehat{U}_{h,4})$, $E_\infty(\widehat{U}_h)$, $E_\infty(u_h)$
-- right, $N_{\sin}=512$}
\label{log10 E N=512}
\end{figure}

We use 4 different coefficients $a_M(\cdot)$ with different intensities of 
oscillation: $N_{\sin}=64$,$128$,$256$,$512$ (see Fig.\ref{aM N=64,128}).
From this series we can observe what happens when $a_M(\cdot)$ becomes 
more and more oscillatory, and guess further behaviour towards more realistic situations.
One test case ($N_{\sin}=256$, $h=1/16$) is illustrated in 
Fig.\ref{A N=256,h=1/16 bw},\ref{comp N=256 h=1/16} (see also \cite{LaptevBelouettar},
where similar results for another $a_M(\cdot)$ were presented). 
The curves $c_1$,$c_2$,$c_3$ are plotted in Fig.\ref{log10 E N=64}--
Fig.\ref{log10 E N=512}.
The contrast of the averaged coefficient is presented in Fig.
\ref{Contrast2}.
%-----------------------------------------------------------------------------------
\begin{figure}[h]
\centerline{\includegraphics[width=1.0\linewidth]
{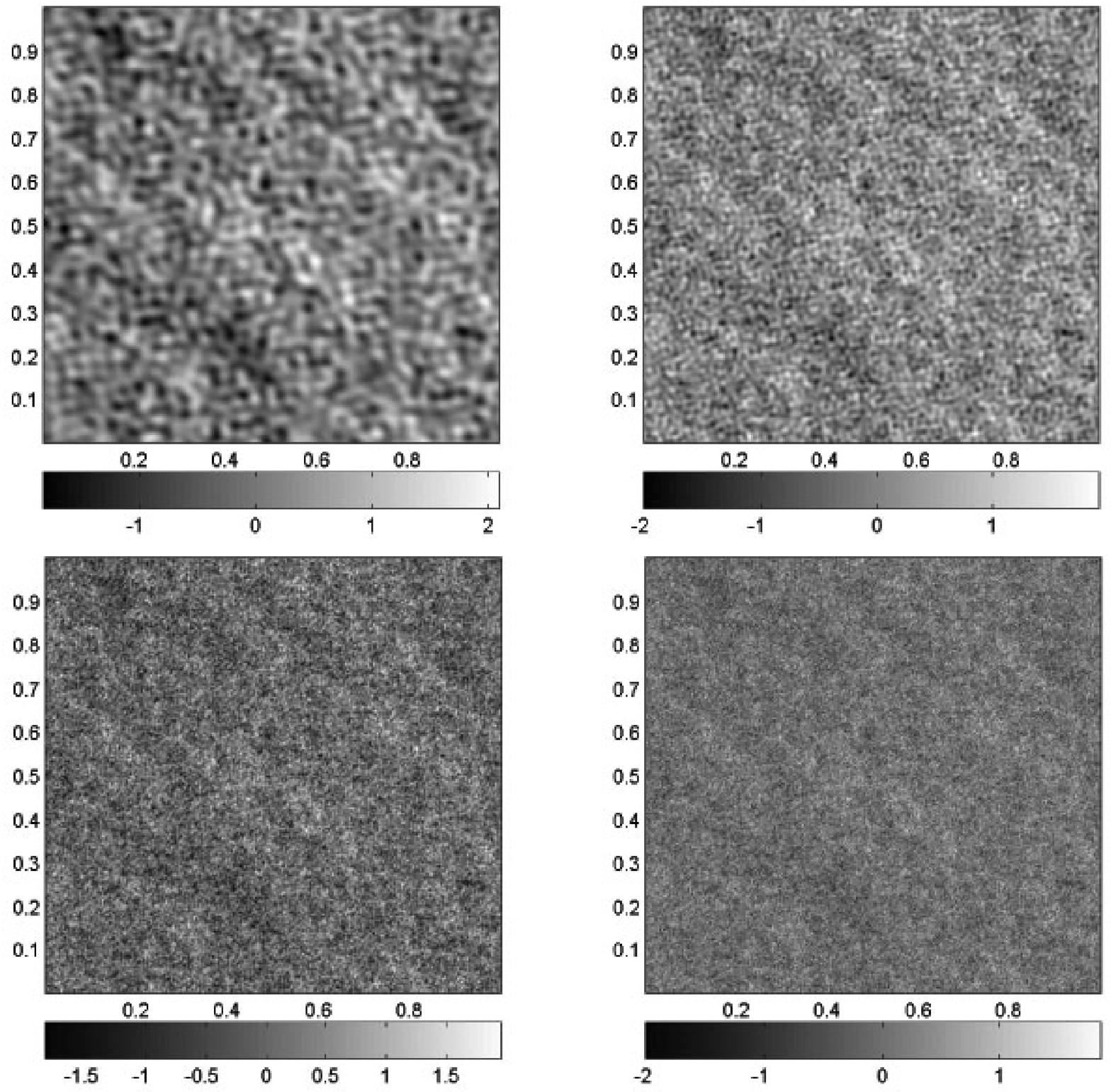}}
\caption{$log_{10}\bigl(a_M(\cdot)\bigr)$ for $N_{\sin}=64$ (top-left) and
  $N_{\sin}=128$ (top-right), $N_{\sin}=256$ (bottom-left) and 
  $N_{\sin}=512$ (bottom-right).}
\label{aM N=64,128}
\end{figure}
%\begin{figure}[h]
%\centerline{\includegraphics[width=1.0\linewidth]
%{./log10a_M64_128_bw.eps}}
%\caption{$log_{10}\bigl(a_M(\cdot)\bigr)$ for $N_{\sin}=64$ (left) and
%  $N_{\sin}=128$ (right).}
%\label{aM N=64,128}
%\end{figure}
%\begin{figure}[h]
%\centerline{\includegraphics[width=1.0\linewidth]
%{./log10a_M256_512_bw.eps}}
%\caption{$log_{10}\bigl(a_M(\cdot)\bigr)$ for $N_{\sin}=256$ (left) and
%  $N_{\sin}=512$ (right).}
%\label{aM N=256,512}
%\end{figure}
%---------------------example RndSinTest256_Big/H16_512_257-----------
\begin{figure}[h]
\centerline{\includegraphics[width=1.0\linewidth]
{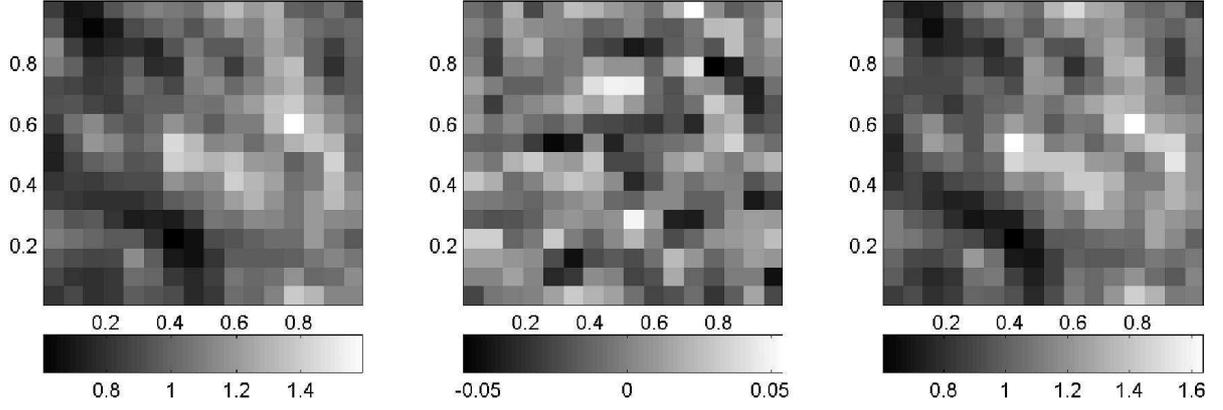}}
\caption{$A_{11}$ (left), $A_{12}=A_{21}$ (middle), $A_{22}$ (right) for $N_{\sin}=256$, $h=1/16$, ${\mathcal D}_2$-extension.}
\label{A N=256,h=1/16 bw}
\end{figure}
%---------------------------------------------------------------------
%-------------------------------------------------------------------------------
\subsubsection
{An interpretation of the presented 2D results for ${\mathcal D}_2$ extensions} 
With the help of the information presented in
Fig.\ref{Ming Yue test}, Fig.\ref{log10 E N=64}--\ref{log10 E N=512}
it is possible to estimate the abilities of the proposed ${\mathcal D}_2$ 
averaging approach 
($c_2$ -- practical, $c_3$ -- theoretical) 
in comparison with the 
direct numerical approach ($c_1$). 

For each $a_M(\cdot)$ we introduce a level $H_{res}$ which
approximately separates the discretization steps $\{h\}$ into two groups:
1) resolving ($h<H_{res}$) and 2) not resolving
($h>H_{res}$) the initial  coefficient $a_M(\cdot)$.
$H_{res}$ is a characteristic value, it is not uniquely defined.
We can choose $H_{res}\simeq 1/(2\cdot 65)$ for the first 2D test 
(Subsection.\ref{subsection MingYue aM}), and $H_{res}\simeq 1/N_{\sin}$
for the rest 4 tests (Subsection.\ref{subsection Random aM}).

When $h<H_{res}$, $c_1(h)$ is a monotone increasing (with a constant rate) 
function of $h$. In the region $h>H_{res}$, $c_1(h)$ is nearly horizontal
since
the direct numerical methods fail to approximate well problems with rapidly oscillated 
coefficients until the coefficients are resolved
(such behaviour is not shown in our figures, except Fig.\ref{log10 E N=512}). 
 
$c_2(h)$, $c_3(h)$ behave in a more complicated way. 
The upscaling is the most effective for coarse grids, $h>H_{upsc}$,
where $c_2(h)$, $c_3(h)$ are monotone increasing (with a constant rate) functions 
of $h$, almost coinsident to each other. 
To illustrate the choice of $H_{upsc}$, we refer to Fig.\ref{Ming Yue test} and Fig.\ref{log10 E N=512},
where $H_{upsc}\simeq 1/32$ and $H_{upsc}\simeq 1/16$ respectively.
   
When $h$ decreases further, $h<H_{upsc}$, the accuracy of the approximation improves but with the  
slowing down rate. The averaging still makes sense, but it is less effective as before.
In all cases except Fig.\ref{Ming Yue test}, $c_2$ reaches a 
local minimum at some $h=H_{acc}$. Further grid refinement in the averaging process gives deterioration 
in the accuracy. Monotone is a desirable property for the 'accuracy vs. discretization size' functions, 
but unfortunately it is unlikely to hold even for $c_3$ curve.
$c_2(h)$ and $c_3(h)$ are almost the same for $h>H_{dev}$ and start to deviate from each other for 
smaller $h$. This happens since the increasing contrast of $A(\cdot)$ (see Fig.\ref{Contrast2}) 
prevents the accurate solving of (\ref{averaged problem}) on the $h$-grid.

We observe that in the region of the resolved $a_M(\cdot)$,  
$c_2$ comes close to $c_1$ (with similar slope) and possibly crosses it. 
For small enough $h$ ($h=\bar\varepsilon/2< H_{res}$) and continuous $a_M(\cdot)$,
the coefficient used in cell problems has a small variation. Consequently the
averaged coefficient $A(\cdot)$ can be seen as a perturbation of $a_M(\cdot)$.  
Thus, there is no surprise that (\ref{initial problem}), 
(\ref{averaged problem}) after solving on the same $h$-grid by the same numerical 
method lead to similar results for $h <H_{res}$.
Also, we note that it is intuitively better to apply a numerical method directly to 
$a_M(\cdot)$ than to its perturbation $A(\cdot)$ when the grid easily resolves the 
initial coefficient. 
This gives some explanation 
why the averaging algorithms rapidly improving at coarse $h$ have to slow
down and to 'wait' the direct method.
Similar behaviour is called ''resonance'' in the terminology of the multiscale 
finite element method \cite{MFEM}.

Let us look how the curves change when $a_M(\cdot)$ becomes more and more oscillatory
($N_{\sin}$ increases from $64$ in Fig.\ref{log10 E N=64} to $512$ in Fig.\ref{log10 E N=512}):
1) $c_1$ moves to the left -- $H_{res}$ decreases; 
2) the region where the averaging is effective has a tendency to expand -- 
$H_{upsc}$, $c_2(H_{upsc})$ decrease;
3) improving of the best accuracy which can be achived on coarse grids 
(it can be roughly characterized by $c_2(H_{acc})$ if the local minimum exists).
 
The quantity
$$
C_A=\frac{\sup\limits_{x\in\Omega}\max\{A_{11}(x),A_{22}(x)\}}
{\inf\limits_{x\in\Omega}\min\{A_{11}(x),A_{22}(x)\}}
$$
plotted in Fig.\ref{Contrast2} for different $\bar\varepsilon$ and $N_{\sin}$ is related to the 
contrast of $A(\cdot)$. The averaged coefficient $A(\cdot)$ is rapidly oscillated
when $\bar\varepsilon$ is small, and $A(\cdot)\simeq const$ when $\bar\varepsilon$ is large.
In other words, $A(x_1)\simeq A(x_2)$ even if $W_{x_1}\cap W_{x_2}=\emptyset$
and $x_1$ and $x_2$ are far from each other. 
This could be an indication of some statistical properties
of our coefficients $a_M(\cdot)$, possibly useful for reducing the computational cost of the 
averaging (see the discussion of linear and sub-linear cost of upscaling algorithms in 
\cite{EE},\cite{MingYue}).
%------------------------------------------------------------------------------
\section{Conclusion}
In this article the averaging algorithm for the second order elliptic equation
using ${\mathcal C}$ and ${\mathcal D}_k$ two-scale extensions was described 
in details and applied to several one and two dimensional model problems.
Our purpose was to show that there are non-periodic coefficients $a_M(\cdot)$
for which the standard periodic homogenization together with the 
two-scale extensions could provide reasonably good averaged coefficients.
For the test cases we investigated how the quality of the approximation
depends on the averaging size $\bar\varepsilon$, and how the averaged approximations
$U_h$ and $\widehat{U}_h$ perform against
the direct numerical approximation (without averaging) $u_h$.

We need to mention that one can construct such initial coefficients 
$a_M(\cdot)$ for which the presented here averaging algorithm fails to approximate well 
on coarse grids. In these cases the averaging has no advantage over the direct 
numerical method. The topic we are planning to address in a forthcoming work.

\end{document}